\documentclass[isolatin1]{article}
\usepackage{amsmath,amsthm,mathtools,mathrsfs,amsfonts,amssymb,stmaryrd,xy}
\usepackage{mathabx}  
\usepackage[a4paper]{geometry}
\theoremstyle{definition}
\newtheorem{remark}{Remark}[section]

\newtheorem{example}[remark]{Example}
 \theoremstyle{plain}
\newtheorem{definition}[remark]{Definition}

\newtheorem{proposition}[remark]{Proposition}

\newtheorem{lemma}[remark]{Lemma}

\newtagform{empty}{}{}
\newcommand{\complex}{\mathbb{C}}

\newcommand{\naturals}{\mathbb{N}}

\DeclareMathOperator{\Id}{id}

\DeclareMathOperator{\Ad}{Ad}

\DeclareMathOperator{\Mor}{Mor}

\newcommand{\lnspan}{\big[} \newcommand{\rnspan}{\big]}

\newcommand{\frakB}{\mathfrak{B}}
\newcommand{\frakBo}{\frakB^{\dag}}
\newcommand{\frakC}{\mathfrak{C}}
\newcommand{\frakCo}{\mathfrak{C}^{\dag}}

\newcommand{\frakH}{\mathfrak{H}}
\newcommand{\frakK}{\mathfrak{K}}

\newcommand{\cbasel}[2]{(\mathfrak{#2}, \mathfrak{#1}, \mathfrak{#1}^{\dag})}
\newcommand{\cbases}[2]{{_{\mathfrak{#1}}}\mathfrak{#2}_{\mathfrak{#1}^{\dag\!}}} 
\newcommand{\cbaseos}[2]{{_{\mathfrak{#1}^{\dag\!}}}\mathfrak{#2}_{\mathfrak{#1}}} 
\newcommand{\cbasesb}{\cbases{B}{H}}
\newcommand{\cbasesc}{\cbases{C}{K}}

\newcommand{\cbaseosb}{\cbaseos{B}{H}}

\newcommand{\lt}{\smalltriangleleft}
\newcommand{\rt}{\smalltriangleright}

\newcommand{\hbeta}{\widehat{\beta}}
\newcommand{\hsigma}{\widehat{\sigma}}

\newcommand{\frei}{\,\cdot\,}
\newcommand{\mycong}{\xrightarrow{\cong}}

\newcommand{\rtensor}[3]{ {_{#1}\! \underset{#2}{\otimes}\! {}_{#3}}}

\newcommand{\htensor}[2]{\rtensor{#1}{\frakH}{#2}}

\newcommand{\mtimes}{\underset{\mu}{\otimes}}
\newcommand{\motimes}{\underset{\mu^{op}}{\otimes}}

\newcommand{\rtensorab}{\htensor{\alpha}{\beta}}

\newcommand{\rtensorrs}{\rtensor{\rho}{\mu}{\sigma}}

\newcommand{\rtensorh}{\underset{\frakH}{\otimes}}
\newcommand{\rtensorm}{\underset{\mu}{\otimes}}

\newcommand{\fibre}[3]{ {_{#1}\! \underset{#2}{\ast}\! {}_{#3}}}

\newcommand{\rfibrers}{\fibre{\rho}{\mu}{\sigma}}

\newcommand{\fsource}{\htensor{\hbeta}{\alpha}}
\newcommand{\frange}{\rtensorab}

\newcommand{\fibreab}{\fibre{\alpha}{\frakH}{\beta}}
\newcommand{\sfsource}{\fsource}
\newcommand{\sfrange}{\frange}

\newcommand{\tl}{\ensuremath \olessthan}
\newcommand{\tr}{\ensuremath \ogreaterthan}

\newcommand{\Hsource}{H \fsource H}
\newcommand{\Hrange}{H \frange H}

\newcommand{\pHsource}{H\rtensor{\hsigma}{\mu}{\rho} H}
\newcommand{\pHrange}{H\rtensor{\rho}{\mu^{op}}{\sigma} H}

\newcommand{\sHsource}{H \sfsource H}
\newcommand{\sHrange}{H \sfrange H}

\newcommand{\HfibreK}{H \rtensorab K}

\newcommand{\Hone}{H \sfsource H \sfsource H}

\newcommand{\Htwo}{H \sfrange H \sfsource  H}

\newcommand{\Hthree}{H \sfrange H \sfrange H}

\newcommand{\Hfour}{(\sHsource) \rtensor{\alpha \lt \alpha}{\frakH}{\beta} H}

\newcommand{\Hfive}{H \rtensor{\hbeta}{\frakH}{\alpha \rt \alpha} (\sHrange)}

\newcommand{\Hfourlt}{H \sfsource H \rtensor{\beta}{\frakH}{\alpha} H}
\newcommand{\Hfourrt}{\big(\sHrange\big) \rtensor{\hbeta \lt
    \beta}{\frakH}{\alpha} H}

\newcommand{\AfibreB}{A \fibreab B}

\newcommand{\kalpha}[1]{|\alpha\rangle_{\leg{#1}}}
\newcommand{\balpha}[1]{\langle\alpha|_{\leg{#1}}}
\newcommand{\kbeta}[1]{|\beta{}\rangle_{\leg{#1}}}
\newcommand{\bbeta}[1]{\langle\beta|_{\leg{#1}}}

\newcommand{\cfact}{\ensuremath \mathrm{C}^{*}\mathrm{\text{-}fact}}

\newcommand{\leg}[1]{[#1]}

\title{Finite-dimensional Hopf $C^{*}$-bimodules and
  $C^{*}$-pseudo-multiplicative unitaries}

\author{Thomas Timmermann\\[1ex]
  \texttt{timmermt@math.uni-muenster.de}\\ SFB 478
  ``Geometrische Strukturen in der Mathematik''\\ Hittorfstr.\
  27, 48149 M\"unster}

\date{\today}

\begin{document}
 \xyrequire{matrix} \xyrequire{arrow}

\maketitle 

\abstract{Finite quantum groupoids can be described in many equivalent
  ways \cite{nikshych:algversions,schauenburg,vallin:3}: In terms of the weak Hopf
  $C^{*}$-algebras of B\"ohm, Nill, and Szlach\'anyi \cite{boehm0} or the
  finite-dimensional Hopf-von Neumann bimodules of Vallin
  \cite{vallin:1}, and in terms of finite-dimensional multiplicative
  partial isometries \cite{boehm:mpi} or the finite-dimensional
  pseudo-multiplicative unitaries of Vallin \cite{vallin:2}.

  In this note, we show that in finite dimensions, the notions of a
  Hopf-von Neumann bimodule and of a pseudo-multiplicative unitary
  coincide with the notions of a concrete Hopf-$C^{*}$-bimodule and of
  a $C^{*}$-pseudo-multiplicative unitary, respectively.
}

\section{Introduction}

The theory of quantum groupoids is very well understood in the finite
and in the measurable case, that is, in the setting of finite-dimensional
$C^{*}$-algebras
\cite{vainer,schauenburg,vallin:3} and in the setting
of von Neumann algebras \cite{lesieur}. The basic objects in this
theory are the weak Hopf $C^{*}$-algebras and the multiplicative
partial isometries of B\"ohm, Nill, and Szlach\'anyi
\cite{boehm0,boehm1,boehm2,boehm:mpi} on one side and the Hopf-von
Neumann bimodules and the pseudo-multiplicative unitaries of Vallin
\cite{vallin:2} on the other side. For finite quantum groupoids, both
approaches are well known to be equivalent
\cite{nikshych:algversions,vallin:3}.

To extend the theory of quantum groupoids to the locally compact case,
that is, to the setting of $C^{*}$-algebras, we introduced the notion
of a concrete Hopf $C^{*}$-bimodule and of a
$C^{*}$-pseudo-multiplicative unitary \cite{timmer:cpmu}. In this
short note, we show that in the finite-dimensional case, these
concepts coincide with the notion of a Hopf-von Neumann bimodule and
of a pseudo-multiplicative unitary, respectively.  This note is of
expository nature and the results contained in it are straightforward.

This work was supported by the SFB 478 ``Geometrische Strukturen in
der Mathematik'' which is funded by the Deutsche
  Forschungsgemeinschaft (DFG).

  \paragraph{Organization} We proceed as follows: 

  In Section 2, we show that every $C^{*}$-factorization of a
  finite-dimensional Hilbert space is uniquely determined by the
  associated representation, and that the $C^{*}$-relative tensor
  product of finite-dimensional Hilbert spaces introduced in
  \cite{timmer:cpmu} coincides with the usual relative tensor product.

  In Section 3, we show that in the finite-dimensional case, the
spatial fiber product of $C^{*}$-algebras introduced in
  \cite{timmer:cpmu} coincides with the usual fiber product of von
  Neumann algebras, and that the notion of a finite-dimensional
  concrete Hopf $C^{*}$-bimodule and of a finite-dimensional Hopf-von
  Neumann bimodule are equivalent. 

  In Section 4, we show that for finite-dimensional Hilbert spaces,
  the notion of a $C^{*}$-pseudo-multiplicative unitary
  \cite{timmer:cpmu} and of a pseudo-multiplicative unitary are
  equivalent, and remark that the associated concrete Hopf
  $C^{*}$-bimodules and Hopf-von Neumann bimodules coincide.

\paragraph{Preliminaries}

Given a subset $Y$ of a normed space $X$, we denote by $[Y] \subseteq X$
the closed linear span of $Y$.

Given a Hilbert space $H$ and a subset $X \subseteq {\cal L}(H)$, we
denote by $X'$  the commutant of $X$. Given  Hilbert spaces $H$, $K$,
a $C^{*}$-subalgebra $A \subseteq {\cal L}(H)$, and a $*$-homomorphism
$\pi \colon A \to {\cal L}(K)$, we put
\begin{align*}
  {\cal L}^{\pi}(H,K) := \{ T \in {\cal L}(H,K) \mid Ta = \pi(a)T
  \text{ for all } a \in A\};
\end{align*}
thus, for example, $A' = {\cal L}^{\Id_{A}}(H)$.

We shall make extensive use of (right) $C^{*}$-modules, also known as
Hilbert $C^{*}$-modules or Hilbert modules. A standard reference is
\cite{lance}.

All sesquilinear maps like inner products of Hilbert spaces
or $C^{*}$-modules are assumed to be conjugate-linear in the first
component and linear in the second one.

Let $A$ and $B$ be $C^{*}$-algebras.  
Given $C^{*}$-modules $E$ and $F$ over $B$, we denote the space of all
adjointable operators $E\to F$ by ${\cal L}_{B}(E,F)$.

Let $E$ and $F$ be $C^{*}$-modules over $A$ and $B$, respectively, and
let $\pi \colon A \to {\cal L}_{B}(F)$ be a $*$-homomorphism. Then one
can form the internal tensor product $E \otimes_{\pi} F$, which is a
$C^{*}$-module over $B$ \cite[Chapter 4]{lance}. This $C^{*}$-module
is the closed linear span of elements $\eta \otimes_{A} \xi$, where
$\eta \in E$ and $\xi \in F$ are arbitrary, and $\langle \eta
\otimes_{\pi} \xi|\eta' \otimes_{\pi} \xi'\rangle = \langle
\xi|\pi(\langle\eta|\eta'\rangle)\xi'\rangle$ and $(\eta \otimes_{\pi}
\xi)b=\eta \otimes_{\pi} \xi b$ for all $\eta,\eta' \in E$, $\xi,\xi'
\in F$, and $b \in B$.  We denote the internal tensor product by
``$\tr$''; thus, for example, $E \tr_{\pi} F=E \otimes_{\pi} F$. If
the representation $\pi$ or both $\pi$ and $A$ are understood, we
write ``$\tr_{A}$'' or  ``$\tr$'', respectively, instead of
$"\tr_{\pi}$''. 

Given $E$, $F$ and $\pi$ as above, we define a {\em flipped internal
  tensor product} $F {_{\pi}\tl} E$ as follows. We equip the algebraic
tensor product $F \odot E$ with the structure maps $\langle \xi \odot
\eta | \xi' \odot \eta'\rangle := \langle \xi| \pi(\langle
\eta|\eta'\rangle) \xi'\rangle$, $(\xi \odot \eta) b := \xi b \odot
\eta$, and by factoring out the null-space of the semi-norm
$\zeta\mapsto \| \langle \zeta|\zeta\rangle\|^{1/2}$ and taking
completion, we obtain a $C^{*}$-$B$-module $F {_{\pi}\tl} E$.  This is
the closed linear span of elements $\xi {_{\pi}\tl} \eta$, where $\eta \in E$
and $\xi \in F$ are arbitrary, and $\langle \xi {_{\pi}\tl} \eta|\xi'
{_{\pi}\tl} \eta'\rangle = \langle \xi|\pi(\langle\eta|\eta'\rangle)\xi'\rangle$ and
$(\xi {_{\pi}\tl} \eta)b=\xi b {_{\pi}\tl} \eta$ for all $\eta,\eta' \in E$, $\xi,\xi'
\in F$, and $b\in B$. As above, we write ``${_{A}\tl}$'' or simply
``$\tl$'' instead of ``${_{\pi}\tl}$'' if the representation $\pi$ or
both  $\pi$ and $A$ are understood, respectively.

Evidently, the usual and the flipped internal tensor product are related by
a unitary map $\Sigma \colon F \tr E \mycong E \tl F$, $\eta \tr \xi
\mapsto \xi \tl \eta$.

Given a state $\mu$ on a finite-dimensional $C^{*}$-algebra $B$, we
denote by $(H_{\mu},\pi_{\mu},\zeta_{\mu})$ a GNS-representation for
$\mu$, by $J_{\mu}\colon H_{\mu} \to H_{\mu}$ the modular conjugation
(an antilinear isometry), and by $\pi_{\mu}^{op}\colon B^{op} \to
{\cal L}(H_{\mu})$ the representation given by $b^{op} \mapsto
J_{\mu}\pi_{\mu}(b)^{*}J_{\mu}$.

\section{The relative tensor product of finite-dimensional
  Hilbert spaces }

In the finite-dimensional case, the $C^{*}$-relative tensor product
and the usual fiber product of Hilbert spaces coincide. Before we can
prove this assertion, we  need to recall the notion of a
$C^{*}$-base and of a $C^{*}$-factorization.

\paragraph{$C^{*}$-bases}
Recall that a {\em $C^{*}$-base} is a triple
$(\frakH,\frakB,\frakBo)$, shortly written $\cbases{B}{H}$, consisting
of a Hilbert space $\frakH$ and two commuting nondegenerate
$C^{*}$-algebras $\frakB,\frakB^{\dag} \subseteq {\cal L}(\frakH)$.
We say that two $C^{*}$-bases $\cbasel{B}{H}$ and $\cbasel{C}{K}$ are
equivalent if there exists a unitary $U \colon H \to K$ such that
$\frakC = \Ad_{U}(\frakB)$ and $\frakCo = \Ad_{U}(\frakBo)$; in that
case, we write $\cbasel{B}{H} \underset{U}{\sim} \cbasel{C}{K}$.
\begin{definition}
  Let $\cbasesb$ be a $C^{*}$-base. We call a vector $\zeta \in
  \frakH$ {\em bicyclic} if it is cyclic for $\frakB$ and for
  $\frakBo$. We call $\cbasesb$ {\em standard} if there exists a
  bicyclic vector $\zeta \in \frakH$, and {\em finite-dimensional} if
  $\frakH$ has finite dimension.
\end{definition}
\begin{example} \label{example:gns-cbase}
  If $\mu$ is a KMS-state on a $C^{*}$-algebra $N$, then the triple
  $(H_{\mu},\pi_{\mu}(N),\pi_{\mu}^{op}(N^{op}))$ is a standard
  $C^{*}$-base, called the {\em $C^{*}$-base associated to $\mu$}, and
  $\zeta_{\mu} \in H_{\mu}$ is bicyclic.
\end{example}
\begin{lemma}
  If $\cbasesb$ is standard and finite-dimensional, then $\frakBo=\frakB'$ and
  $\frakB=(\frakBo)'$. 
\end{lemma}
\begin{proof}
  By definition, $\frakBo \subseteq \frakB'$. If $\zeta \in \frakH$ is
  bicyclic, then the map $j \colon \frakB' \to H$ given by $T \mapsto
  T\zeta$ is injective and $j(\frakBo)=\frakH =
  j(\frakB')$. Therefore $\frakBo = \frakB'$.
\end{proof}

Using standard results on   GNS-representations and the lemma above,
one finds:
\begin{lemma} \label{lemma:cbase-gns}
  If $\cbasesb$ is a finite-dimensional standard $C^{*}$-base and
  $\zeta \in \cbasesb$ is bicyclic, then the state $\mu:=\langle \zeta
  | \frei\zeta\rangle$ on $\frakB$ is faithful and there exists a
  unique unitary $U\colon \frakH \to H_{\mu}$ such that $U\zeta =
  \zeta_{\mu}$ and $\pi_{\mu}(\frakB)=\Ad_{U}(\frakB)$; moreover,  then
  $\pi_{\mu}^{op}(\frakB^{op})=\frakBo$. \qed
\end{lemma}
\begin{remark} \label{remark:frakbo-frakbop} Let $\cbasesb$ be a
  finite-dimensional standard $C^{*}$-base, $\zeta \in \cbasesb$
  bicyclic, and $U \colon \frakH \to H_{\mu}$ as above. Then we can
  identify $\frakBo$ with $\frakB^{op}$ via $(\pi_{\mu}^{op})^{-1}
  \circ \Ad_{U}$. More concretely, if $J\colon \frakH \to \frakH$
  denotes the antiunitary part in the polar decomposition of the map
  $\frakH \to \frakH$, $b\zeta \mapsto b^{*}\zeta$, then
  $J=U^{*}J_{\mu}U$, and the map $b^{op} \mapsto Jb^{*}J$ is an isomorphism
  $\frakB^{op} \mycong \frakBo$.
\end{remark}
\paragraph{$C^{*}$-factorizations}
Let $\cbasesb$ be a standard $C^{*}$-base with bicyclic vector $\zeta
\in \frakH$ and let $H$ be a Hilbert space.

Recall that a {\em $C^{*}$-factorization} of $H$ with respect to
$\cbasesb$ is a closed subspace $\alpha \subseteq {\cal L}(\frakH,H)$
satisfying $[\alpha^{*}\alpha] = \frakB$, $[\alpha \frakB] = \alpha$,
and $[\alpha \frakH] = H$.  We denote by $\cfact(H;\cbasesb)$ the set
of all $C^{*}$-factorizations of $H$ with respect to $\cbasesb$.

\begin{lemma} \label{lemma:dense-embedding} For each $\alpha \in
  \cfact(H;\cbasesb)$, the map $\alpha \to H$, $\xi \mapsto \xi
  \zeta$, is injective and has dense image.
\end{lemma}
\begin{proof}
  If $\xi \in \alpha$ and $\xi \zeta = 0$, then $\xi\frakH = [\xi
  \frakBo \zeta ] = [\rho_{\alpha}(\frakBo)\xi\zeta] = 0$ and hence
  $\xi = 0$. Therefore, the map $\xi \mapsto \xi \zeta$ is
  injective. It has dense image because $[\alpha \zeta] = [\alpha
  \frakB \zeta] = [\alpha \frakH] = H$.
\end{proof}

From now on, we assume that $H$ has finite dimension. Let $\rho \colon
\frakBo \to {\cal L}(H)$ be a nondegenerate faithful representation
and put
\begin{align*}
  {\cal L}^{\rho}(\frakH,H) := \{ T \in {\cal L}(\frakH,H) \mid
  Tb^{\dag} = \rho(b^{\dag})T \text{ for all } b^{\dag} \in \frakBo\}.
\end{align*}
\begin{lemma} \label{lemma:intertwiner}
  \begin{enumerate}
  \item For each $\xi \in H$, there exists a unique
    $R^{\rho}_{\zeta}(\xi) \in {\cal L}^{\rho}(\frakH,H)$ such that
    $R^{\rho}_{\zeta}(\xi)  \zeta = \xi$.
  \item For each $T \in {\cal L}^{\rho}(\frakH,H)$, one has
    $T=R^{\rho}_{\zeta}(T\zeta)$.  
  \end{enumerate}
\end{lemma}
\begin{proof}
  Straightforward.
\end{proof}

The $C^{*}$-factorizations of $H$ are uniquely determined by their
associated representations:
\begin{proposition} \label{proposition:rep-cfact} Let $\rho \colon
  \frakBo \to {\cal L}(H)$ be a nondegenerate faithful representation
  and $\alpha\in \cfact(H;\cbasesb)$.
  \begin{enumerate}
  \item There exists a unique nondegenerate faithful representation
    $\rho_{\alpha} \colon \frakBo \to {\cal L}(H)$ such that
    $\rho_{\alpha}(b^{\dag})\xi\zeta = \xi b^{\dag}\zeta$ for all
    $b^{\dag} \in \frakB^{\dag}$, $\xi\in \alpha$, $\zeta \in \frakH$.
  \item ${\cal L}^{\rho}(\frakH,H) \in \cfact(H;\cbasesb)$.
  \item $\rho_{\alpha}=\rho$ if and only if $\alpha={\cal
      L}^{\rho}(\frakH,H)$.
  \end{enumerate}
\end{proposition}
\begin{proof}
  i) The representation $\rho_{\alpha}$ is well-defined because for
  all $\xi,\xi' \in \alpha$, $\zeta,\zeta' \in \frakH$, and $b^{\dag}
  \in \frakBo$,
  \begin{align*}
    \langle \xi \zeta | \xi' b^{\dag}\zeta'\rangle = \langle \zeta |
    \xi^{*} \xi' b^{\dag}\zeta'\rangle =
    \langle \zeta | b^{\dag}\xi^{*}\xi' \zeta'\rangle =
    \langle \xi(b^{\dag})^{*}\zeta|\xi'\zeta'\rangle;
  \end{align*}
  here, we used $\alpha^{*}\alpha \subseteq \frakB \subseteq
  (\frakBo)'$.  Combining this calculation with the relation $\lnspan
  \alpha^{*}\alpha \frakH\rnspan = \frakH$, we find that
  $\rho_{\alpha}$ is faithful. It is
  nondegenerate because $\lnspan \rho_{\alpha}(\frakBo)H\rnspan =
  \lnspan \alpha \frakBo \frakH\rnspan = \lnspan \alpha \frakH\rnspan
  = H$.

  \smallskip

  ii) Put $\beta:={\cal L}^{\rho}(\frakH,H)$.  Lemma
  \ref{lemma:intertwiner} i) implies $\lnspan \beta \frakH\rnspan =
  H$, and a short calculation shows $\lnspan \beta^{*}\beta\rnspan
  \subseteq (\frakBo)' = \frakB$. We prove that this inclusion is an
  equality. Choose a bicyclic vector $\zeta \in \frakH$ and consider
  the map $j \colon \frakB \to {\cal L}(\frakBo,\complex)$ given by $c
  \mapsto \langle \zeta|c\frei\zeta\rangle$. Since $\zeta$ is cyclic
  for $\frakBo$ and $\frakBo$ commutes with $\frakB$, this map is
  injective. Moreover, since $\rho$ is faithful and
  \begin{align*}
    j(R_{\zeta}^{\rho}(\xi)^{*}R_{\zeta}^{\rho}(\xi')) = \langle
    \zeta|R_{\zeta}^{\rho}(\xi)^{*}R_{\zeta}^{\rho}(\xi')(\frei)
    \zeta\rangle = \langle \xi|\rho(\frei) \xi'\rangle \quad \text{for
      all } \xi,\xi' \in H,
  \end{align*}
  we have $j([\beta^{*}\beta]) = {\cal L}(\frakBo,\complex) \supseteq
  j(\frakB)$. Consequently, $\lnspan\beta^{*}\beta\rnspan=\frakB$.
  Finally, we prove $[\beta \frakB]=\beta$.   Short calculations show
  that $\lnspan \beta\beta^{*}\rnspan \subseteq \rho(\frakBo)'$ and $TR_{\zeta}^{\rho}(\xi)=R_{\zeta}^{\rho}(T\xi)$ for each $T \in
  \rho(\frakBo)'$, $\xi \in H$, and therefore, $[\beta
  \frakB]=[\beta \beta^{*}\beta] \subseteq [\rho(\frakBo)'\beta]=\beta$.
  Conversely,  $\beta=\beta\Id_{\frakH} \subseteq [\beta\frakB]$.

  \smallskip

  iii) If $\alpha={\cal L}^{\rho}(\frakH,H)$, then evidently,
  $\rho_{\alpha}=\rho$. Conversely, assume $\rho_{\alpha}=\rho$. Then
  evidently $\alpha \subseteq {\cal L}^{\rho}(\frakH,H)$. We prove
  that this inclusion is an equality. The map $\alpha \to H$ given by
  $\xi \mapsto \xi \zeta$ is bijective (Lemma
  \ref{lemma:dense-embedding}), so $\dim \alpha = \dim H$. But by
  Lemma \ref{lemma:intertwiner}, $\dim H = \dim {\cal
    L}^{\rho}(\frakH,H)$, and therefore, $\alpha = {\cal
    L}^{\rho}(\frakH,H)$.
\end{proof}

Let $\cbasesc$ be a finite-dimensional standard $C^{*}$-base. Recall
that two $C^{*}$-factorizations $\alpha \in \cfact(H;\cbasesb)$ and
$\beta \in \cfact(H;\cbasesc)$ are {\em compatible} if
$[\rho_{\alpha}(\frakBo)\beta]=\beta$ and
$[\rho_{\beta}(\frakCo)\alpha]=\alpha$. 
\begin{proposition} \label{proposition:rep-fact-compatible} Let $\rho
  \colon \frakBo \to {\cal L}(H)$, $\sigma \colon \frakCo \to {\cal
    L}(H)$ be faithful nondegenerate representations. Then
  $\rho(\frakBo)$ commutes with $\sigma(\frakCo)$ if and only if the
  $C^{*}$-factorizations ${\cal L}^{\rho}(\frakH,H)$ and ${\cal
    L}^{\sigma}(\frakK,H)$ are compatible.
\end{proposition}
\begin{proof}
  Put $\alpha:={\cal L}^{\rho}(\frakH,H)$ and $\beta:={\cal
    L}^{\sigma}(\frakK,H)$; then $\rho(\frakBo) =
  \rho_{\alpha}(\frakBo)$ and $\sigma(\frakCo)=\rho_{\beta}(\frakCo)$
  by Proposition \ref{proposition:rep-cfact}.  If $\alpha$ and $\beta$
  are compatible, then $\rho(\frakBo)$ commutes with $\sigma(\frakCo)$
  by \cite[Remark 2.6]{timmer:cpmu}.  Conversely, if
  $\rho(\frakBo)$ and $\sigma(\frakCo)$ commute, then $\lnspan
  \rho(\frakBo) \beta \rnspan\subseteq \lnspan \sigma(\frakCo)' {\cal
    L}^{\sigma}(\frakK,H)\rnspan = {\cal L}^{\sigma}(\frakK,H)=\beta =
  \Id_{H}\beta \subseteq \lnspan \rho(\frakBo)\beta\rnspan$ and likewise
  $\lnspan \rho_{\beta}(\frakCo) \alpha \rnspan = \alpha$.
\end{proof}

\paragraph{The relative tensor product of Hilbert spaces}
Let us recall the construction of the $C^{*}$-relative tensor product
\cite{timmer:cpmu} and  the usual
relative tensor product of finite-dimensional Hilbert spaces. Suppose that
 \begin{enumerate}
 \item  $H$ and $K$ are finite-dimensional Hilbert spaces,
 \item $N$ is a finite-dimensional $C^{*}$-algebra with a faithful
   state $\mu$ and nondegenerate faithful representations $\rho \colon
   N^{op} \to {\cal L}(H)$, $\sigma \colon N \to {\cal L}(K)$,
 \item $\cbasesb$ is a finite-dimensional standard $C^{*}$-base with
    bicyclic vector $\zeta \in \frakH$ and $C^{*}$-factorizations
    $\alpha \in \cfact(H;\cbasesb)$, $\beta \in
    \cfact(K;\cbaseosb)$,
 \end{enumerate}
such that
\begin{gather} \label{eq:rtp-iso-cond}
  \begin{gathered}
    (\frakH,\frakB,\frakBo) \underset{U}{\sim}
    (H_{\mu},\pi_{\mu}(N),\pi_{\mu}^{op}(N^{op})), \qquad U\zeta =
    \zeta_{\mu},  \\
    \rho = \rho_{\alpha} \circ \Ad_{U^{*}} \circ \pi_{\mu}^{op}, \qquad \sigma
    = \rho_{\beta} \circ \Ad_{U^{*}} \circ \pi_{\mu}.
  \end{gathered}
  \end{gather}
Note that by Example \ref{example:gns-cbase}, Lemma
\ref{lemma:cbase-gns}, and Proposition \ref{proposition:rep-cfact},
given the data listed in ii), we can
construct the data listed in iii) such that
\eqref{eq:rtp-iso-cond} is satisfied, and vice versa.

The relative tensor product of $H$ and $K$ with respect to
$\mu,\rho,\sigma$ is defined as follows.  For each $\xi \in H$ and
$\eta \in K$, there exist unique operators
\begin{align*}
  R_{\mu^{op}}^{\rho}(\xi) &\colon H_{\mu} \to H, \
  \pi_{\mu}^{op}(b^{op})\zeta_{\mu} \mapsto \rho(b^{op})\xi,
  &&\text{and} & R_{\mu}^{\sigma}(\eta) &\colon H_{\mu} \to K, \
  \pi_{\mu}(b)\zeta_{\mu} \mapsto \sigma(b)\eta.
\end{align*}
 Define a sesquilinear form $\langle
\frei | \frei\rangle_{\mu}$ on $H \otimes K$ by
\begin{align*}
  \langle \xi \otimes \eta|\xi' \otimes \eta'\rangle_{\mu} :=
  \big\langle
  R_{\mu^{op}}^{\rho}(\xi')^{*}R_{\mu^{op}}^{\rho}(\xi)\zeta_{\mu}\big|
  R_{\mu}^{\sigma}(\eta)^{*}R_{\mu}^{\sigma}(\eta')\zeta_{\mu}\big\rangle
  \quad \text{for all } \xi,\xi' \in H,\, \eta,\eta' \in K.
\end{align*}
Factoring out the null space of the associated seminorm, we obtain a
Hilbert space $H \rtensorrs K$. For all $\xi \in H$
and $\eta \in K$, we denote by $\xi \rtensor{\rho}{\mu}{\sigma} \eta$
the image of $\xi \otimes \eta$ in $H \rtensorrs K$.

The $C^{*}$-relative tensor product of $H$ and $K$ with respect to
$\cbasesb, \alpha, \beta$ is the internal tensor product $H
\rtensorab K := \alpha \tr \frakH \tl \beta$ \cite[Section
2]{timmer:cpmu}.

\begin{proposition} \label{proposition:rtp-iso} There exists a 
  unitary
  \begin{align*}
    \Phi^{U,\zeta}_{\alpha,\beta} \colon H \rtensorrs K \to H
    \rtensorab K, \quad \xi \rtensorrs \eta \mapsto R^{\rho}_{\mu^{op}}(\xi)U \tr
    \zeta_{\mu} \tl R^{\sigma}_{\mu}(\eta)U.
  \end{align*}
\end{proposition}
\begin{proof}[Proof of Proposition \ref{proposition:rtp-iso}]
   By Proposition \ref{proposition:rep-cfact} and Lemma
   \ref{lemma:intertwiner} ii), $\alpha = \{ R^{\rho}_{\mu^{op}}(\xi)U \mid
  \xi \in H\}$ and $\beta=\{R^{\sigma}_{\mu}(\eta)U \mid \eta \in K\}$, and by
  definition,
  \begin{multline*}
    \langle R^{\rho}_{\mu^{op}}(\xi)U \tr \zeta \tl
    R^{\sigma}_{\mu}(\eta)U | R^{\rho}_{\mu^{op}}(\xi')U \tr \zeta \tl
    R^{\sigma}_{\mu}(\eta')U\rangle = \\ \big\langle \zeta_{\mu}\big|
    R^{\rho}_{\mu^{op}}(\xi)^{*}R^{\rho}_{\mu^{op}}(\xi')
    R^{\sigma}_{\mu}(\eta)^{*}R^{\sigma}_{\mu}(\eta')\zeta_{\mu}\big\rangle
     = \langle \xi \rtensorrs \eta|\xi' \rtensorrs \eta'\rangle
  \end{multline*}
  for all $\xi,\xi'\in H$, $\eta,\eta' \in K$. Therefore, $
  \Phi^{U,\zeta}_{\alpha,\beta}$ is a well-defined isometry. It is
  surjective because $\zeta$ is cyclic for $\frakB$ (and for
  $\frakBo$).
\end{proof}

\section{Finite-dimensional Hopf $C^{*}$-bimodules}

In the finite-dimensional case, the notion of a Hopf $C^{*}$-bimodule
and of a Hopf-von Neumann bimodule coincide. To prove this assertion,
we first review the fiber product of finite-dimensional
$C^{*}$-algebras and the fiber product of morphisms.

\paragraph{The fiber product  of finite-dimensional
  $C^{*}$-algebras}
The spatial fiber product of finite-dimensional $C^{*}$-algebras
\cite{timmer:cpmu} coincides with the usual fiber product of
$C^{*}$-algebras.  To make this statement precise, we briefly recall
the two constructions. Let
\begin{enumerate}
\item $H$ and $K$ be finite-dimensional Hilbert spaces,
\item $A \subseteq {\cal L}(H)$ and $B \subseteq {\cal L}(K)$ be
  nondegenerate $C^{*}$-subalgebras,
  \item $N$ be a finite-dimensional $C^{*}$-algebra with a faithful
    state $\mu$ and injective unital $*$-homomorphisms $\rho
    \colon N^{op} \to A$, $\sigma
    \colon N \to B$,
  \item $\cbasesb$ be a finite-dimensional standard $C^{*}$-base with
    bicyclic vector $\zeta \in \frakH$ and $C^{*}$-factorizations
    $\alpha \in \cfact(A;\cbasesb)$, $\beta \in
    \cfact(B;\cbaseosb)$, where
  \begin{align*}
    \cfact(A;\cbasesb) &= \{ \alpha \in \cfact(H;\cbasesb) \mid
    \rho_{\alpha}(\frakBo) \subseteq A\},  \\
    \cfact(B;\cbaseosb) &= \{ \beta \in \cfact(K;\cbaseosb) \mid
    \rho_{\beta}(\frakB) \subseteq B\},
  \end{align*}
\end{enumerate}
  such that \eqref{eq:rtp-iso-cond} holds.

  The fiber product of $A$ and $B$ with respect to $\mu,\rho,\sigma$
  is defined as follows. For each $S \in A'\subseteq \rho(N^{op})'$
  and $T \in B' \subseteq \sigma(N)'$, there exists a well-defined
  operator
  \begin{align*}
  S \rtensorrs T \colon H \rtensorrs K \to H \rtensorrs K, \quad\xi
  \rtensorrs \eta \mapsto S\xi \rtensorrs T\eta.
\end{align*}
The fiber product of $A$ and $B$ is the commutant $A \rfibrers B :=
(A' \rtensorrs B')' \subseteq {\cal L}(H \rtensorrs K)$.

The spatial $C^{*}$-fiber product of $A$ and $B$ with respect to
$\cbasesb,\alpha,\beta$ is defined as follows.  Using the isomorphisms
  \begin{gather} \label{eq:rtp-space} \alpha \tr_{\rho_{\beta}} K
    \cong \HfibreK \cong H {_{\rho_{\alpha}} \tl} \beta, \quad \xi \tr
    \eta \zeta \equiv \xi \tr \zeta \tl \eta \equiv \xi \zeta \tl
    \eta,
\end{gather}
(see \cite[Section 2]{timmer:cpmu}), one defines for each $\xi \in
\alpha$ and $\eta \in \beta$ operators
\begin{gather} \label{eq:ketbra-ops}
  \begin{aligned}
    |\xi\rangle_{\leg{1}} \colon K &\to \HfibreK, \ \zeta \mapsto \xi
    \tr \zeta, & \langle
    \xi|_{\leg{1}}:=|\xi\rangle_{\leg{1}}^{*}\colon \xi' \tr \zeta
    &\mapsto
    \rho_{\beta}(\langle\xi|\xi'\rangle)\zeta, \\
    |\eta\rangle_{\leg{2}} \colon H &\to \HfibreK, \ \zeta \mapsto
    \zeta \tl \eta, & \langle\eta|_{\leg{2}} :=
    |\eta\rangle_{\leg{2}}^{*} \colon \zeta \tl\eta &\mapsto
    \rho_{\alpha}(\langle \eta|\eta'\rangle)\zeta.
  \end{aligned}
\end{gather}
Put $\kalpha{1} := \big\{
|\xi\rangle_{\leg{1}} \,\big|\, \xi \in \alpha\big\}$ and similarly
define $\balpha{1}$, $\kbeta{2}$, $\bbeta{2}$. Then 
\begin{align*}
  A \fibreab B = \big\{ T \in {\cal L}(H \rtensorab K) \, \big|\, T\kalpha{1},
  T^{*}\kalpha{1} \subseteq [\kalpha{1}B]  \text{ and } T\kbeta{2},
  T^{*}\kbeta{2} \subseteq [\kbeta{2}A]\big\}.
\end{align*}

The two constructions described above coincide in the following sense:
\begin{proposition} \label{proposition:fp-iso} Conjugation by
  $\Phi^{U,\zeta}_{\alpha,\beta} \colon H \rtensorrs K \to H
  \rtensorab K$ induces an isomorphism
  \begin{align*}
    \phi^{U,\zeta}_{\alpha,\beta} \colon A \rfibrers B \mycong A \fibreab B.
  \end{align*}
\end{proposition}
\begin{proof} Put $\Phi:=\Phi^{U,\zeta}_{\alpha,\beta}$ and let $T \in
  {\cal L}(H \rtensorab K)$. By definition, $T \in \Ad_{\Phi}(A
  \rfibrers B)$ if and only if $[T, A' \tl \Id_{\beta}] = 0 = [T,
  \Id_{\alpha} \tr B']$, that is, if and only if for all $\eta,\eta'
  \in \beta$ and $\xi,\xi' \in \alpha$,
  \begin{align*}
    \langle \eta|_{\leg{2}} [T,A' \tl
    \Id_{\beta}]|\eta'\rangle_{\leg{2}} = 0 \quad \text{and} \quad \langle
    \xi|_{\leg{1}} [T,\Id_{\alpha} \tr B']|\xi'\rangle_{\leg{1}} = 0,
  \end{align*}
  or, equivalently, if and only if $\bbeta{2}T\kbeta{2} \subseteq
   A''=A$ and $\balpha{1}T\kalpha{1} \subseteq B''=B$. 

  If $T \in A \fibreab B$, then  $T \in \Ad_{\Phi}(A
  \rfibrers B)$ because
  \begin{align*}
  \bbeta{2}T\kbeta{2} \subseteq
  [\bbeta{2}\kbeta{2}A] = [\rho_{\alpha}(\frakBo)A]=A  
  \end{align*}
  and similarly $\balpha{1}T\kalpha{1} \subseteq B$.

Conversely, if $T \in \Ad_{\Phi}(A \rfibrers B)$, then  $T
  \in A \fibreab B$. Indeed, then
  $\Id_{\alpha} \in [\alpha\alpha^{*}]$ implies
  \begin{align*}
    T\kalpha{1} \in [\kalpha{1} \balpha{1} T \kalpha{1}] \subseteq
    [\kalpha{1} B],
  \end{align*}
  and similarly $T^{*}\kalpha{1} \subseteq [\kalpha{1}B]$ and
  $T\kbeta{2}, T^{*}\kbeta{2} \subseteq [\kbeta{2}A]$.
\end{proof}

\paragraph{Morphisms of finite-dimensional $C^{*}$-algebras}
Let $\cbasesb$ be a $C^{*}$-base.  A {\em nondegenerate
  finite-dimensional concrete (shortly nfc.)
  $C^{*}$-$\cbasesb$-algebra} $(H,A,\alpha)$ consists of a
finite-dimensional Hilbert space $H$, a nondegenerate $C^{*}$-algebra
$A \subseteq {\cal L}(H)$, and a $C^{*}$-factorization $\alpha \in
\cfact(A;\cbasesb)$.  A {\em morphism} of nfc.\
$C^{*}$-$\cbasesb$-algebras $(H,A,\alpha)$ and $(K,B,\beta)$ is a
$*$-homomorphism $\pi \colon A \to B$ such that $\beta =
[I_{\pi}\alpha]$, where $I_{\pi} := \big\{ V \in {\cal L}^{\pi}(H,K)
\,\big|\ V\alpha \subseteq \beta,\, V^{*}\beta \subseteq \alpha\big\}$ \cite{timmer:cpmu,timmer:ckac}. We denote the set of
such morphisms by $\Mor(A_{\alpha},B_{\beta})$.

\begin{lemma} \label{lemma:morphism}
  Let $\cbasesb$ be a standard $C^{*}$-base, $(H,A,\alpha)$,
  $(K,B,\beta)$ nfc.\ $C^{*}$-$\cbasesb$-algebras, and $\pi \colon
  A\to B$ a unital $*$-homomorphism. Then $\pi \in
  \Mor(A_{\alpha},B_{\beta})$ if and only if
  \begin{align} \label{eq:morphism-cond}
   \pi(\rho_{\alpha}(b^{\dag}))=\rho_{\beta}(b^{\dag}) \quad \text{for
     all }
   b^{\dag}\in \frakBo. 
  \end{align}
\end{lemma}
\begin{proof}
  If $\pi$ is a morphism, then \eqref{eq:morphism-cond} holds by
  \cite[Lemma 2.2]{timmer:ckac}. Conversely, assume
  \eqref{eq:morphism-cond}. By Lemma \ref{lemma:intertwiner}, it
  suffices to prove $\lnspan {\cal L}^{\pi}(H,K) H\rnspan = K$. Since
  $A$ has finite dimension, there exist $n \in \naturals$ and central projections
  $p_{1},\ldots,p_{n} \in A$ such that $\sum_{i} p_{i}=1_{A}$ and such
  that each $p_{i}A$ is a matrix algebra.  Since $\sum_{i} \pi(p_{i})
  K=K$, it suffices to show that $\lnspan {\cal
    L}^{\pi_{i}}(p_{i}H,\pi(p_{i})K)p_{i}H\rnspan = \pi(p_{i})K$ for
  each $i=1,\ldots,n$; here, $\pi_{i} \colon p_{i}A \to {\cal
    L}(\pi(p_{i})K)$ denotes the restriction of $\pi$. But both the
  identity representation and the representation $\pi_{i}$ of $p_{i}A$
  are direct sums of the  irreducible representation of the matrix
  algebra $p_{i}A$ which is unique up to unitary equivalence, and
  therefore, $\lnspan {\cal
    L}^{\pi_{i}}(p_{i}H,\pi(p_{i})K)p_{i}H\rnspan = \pi(p_{i})K$.
\end{proof}

\paragraph{The fiber product of morphisms}
 In the finite-dimensional case, the classical fiber
product coincides with the spatial fiber product also on the level of
morphisms. More precisely, let
\begin{enumerate}
\item $H$, $K$, $L$, $M$ be finite-dimensional Hilbert spaces,
\item $A \subseteq {\cal L}(H)$, $B \subseteq {\cal L}(K)$, $C
  \subseteq {\cal L}(L)$, $D \subseteq {\cal L}(M)$ be nondegenerate
  $C^{*}$-algebras,
\item $\phi \colon A \to C$ and $\psi \colon B \to D$ be
  unital $*$-homomorphisms,
\item $N$ be a $C^{*}$-algebra, $\mu$ a faithful state on $N$, and
  $\rho \colon N^{op} \to A$, $\sigma \colon N \to B$, $\upsilon \colon
  N^{op} \to C$, $\omega \colon N \to D$ injective unital
  $*$-homomorphisms,
\item $\cbasesb$ be a standard $C^{*}$-base, $\zeta \in \frakH$ a
  bicyclic vector, and $\alpha \in \cfact(A;\cbasesb)$, $\beta \in
  \cfact(B;\cbaseosb)$, $\gamma \in \cfact(C;\cbasesb)$, $\delta \in
  \cfact(D;\cbaseosb)$,
\end{enumerate}
and assume 
\begin{gather} \label{eq:fp-iso-cond}
  \begin{gathered}
    (\frakH,\frakB,\frakBo) \underset{U}{\sim}
    (H_{\mu},\pi_{\mu}(N),\pi_{\mu}^{op}(N^{op})), \qquad U\zeta =
    \zeta_{\mu},
    \\
    \rho = \rho_{\alpha} \circ \Ad_{U^{*}} \circ \pi_{\mu}^{op}, \quad
    \sigma =\rho_{\beta} \circ \Ad_{U^{*}} \circ \pi_{\mu}, \quad
    \upsilon = \rho_{\gamma} \circ \Ad_{U^{*}} \circ \pi_{\mu}^{op},
    \quad \omega = \rho_{\delta} \circ \Ad_{U^{*}} \circ \pi_{\mu}.
  \end{gathered}
\end{gather}
Note that by Example \ref{example:gns-cbase}, Lemma
\ref{lemma:cbase-gns}, and Proposition \ref{proposition:rep-cfact},
given the data listed in iv), we can construct the data listed in v)
such that \eqref{eq:fp-iso-cond} is satisfied, and vice versa.

By Lemma  \ref{lemma:morphism}, the following
 conditions are equivalent:
 \begin{align} \label{eq:morphism-equivariant}
   i) \ \phi \circ \rho = \upsilon,\ \psi \circ \sigma = \omega, \qquad
   ii) \ \phi \in \Mor(A_{\alpha},C_{\gamma}),\ \psi \in
   \Mor(B_{\beta},D_{\delta}).
 \end{align}
 Assume that these conditions hold. Then by \cite[Proof of
 1.2.4]{sauvageot:2} and \cite[Proposition 3.13]{timmer:cpmu},
 respectively, there exist unique $*$-homomorphisms
\begin{align*}
  \phi \fibre{}{\mu}{} \psi \colon A
\fibre{\rho}{\mu}{\sigma} B \to C \fibre{\upsilon}{\mu}{\omega} D
 \quad \text{and} \quad
  \phi \fibre{}{\frakH}{} \psi \colon \AfibreB \to C
  \fibre{\gamma}{\frakH}{\delta} D 
\end{align*}
 such that for all $X \in {\cal
  L}^{\phi}(H,L)$, $Y \in {\cal 
    L}^{\psi}(K,M)$, $S \in A \rfibrers B$,  $T \in \AfibreB$,
\begin{align*}
  (\phi \fibre{}{\mu}{} \psi)(S) \cdot (X \rtensor{}{\mu}{} Y) = (X
  \rtensor{}{\mu}{} Y) \cdot S
 \quad \text{and} \quad
  (\phi \fibre{}{\frakH}{} \psi)(T) \cdot (X \rtensorh Y) = (X
  \rtensorh Y) \cdot T.
\end{align*}
\begin{proposition}
  If condition \eqref{eq:morphism-equivariant} holds, then the following diagram
  commutes:
\begin{align*} 
  \xymatrix@R=15pt{ {A \fibre{\rho}{\mu}{\sigma} B} \ar[r]^{\phi
      \underset{\mu}{\ast} \psi}
    \ar[d]^{\cong}_{\phi^{U,\zeta}_{\alpha,\beta}} & {C
      \fibre{\upsilon}{\mu}{\omega} D}
    \ar[d]_{\cong}^{\phi^{U,\zeta}_{\gamma,\delta}} \\ 
    {\AfibreB} \ar[r]^{\phi \underset{\frakH}{\ast} \psi} & {C
      \fibre{\gamma}{\frakH}{\delta} D.} }
\end{align*}
\end{proposition}
\begin{proof}
  This follows from the definition of $\phi \underset{\mu}{\ast} \psi$
  and $\phi \underset{\frakH}{\ast} \psi$ and the fact that
  $\Phi^{U,\zeta}_{\gamma,\delta} (X \rtensor{}{\mu}{} Y) = (X
    \rtensor{}{\frakH}{} Y) \Phi^{U,\zeta}_{\alpha,\beta}$
    for all $X \in {\cal L}^{\phi}(H,L)$, $Y \in {\cal
      L}^{\psi}(K,M)$.
\end{proof}

\paragraph{Finite-dimensional concrete Hopf $C^{*}$-bimodules and
  Hopf-von Neumann bimodules}
Let us briefly recall  the definition of a concrete Hopf
$C^{*}$-bimodule and of a concrete Hopf-von Neumann
bimodule.  Suppose that
  \begin{enumerate}
  \item $H$ is a finite-dimensional Hilbert space, $A \subseteq
    {\cal L}(H)$ is a nondegenerate $C^{*}$-subalgebra,
  \item $N$ is a finite-dimensional $C^{*}$-algebra with a faithful
    state $\mu$ and injective unital $*$-homomorphisms $\rho \colon
    N^{op} \to A$ and $\sigma \colon N \to A$ such that $\rho(N^{op})$
    and $\sigma(N)$ commute,
  \item $\cbasesb$ is a finite-dimensional standard $C^{*}$-base with
    bicyclic vector $\zeta \in \frakH$ and compatible $C^{*}$-factorizations
    $\alpha \in \cfact(A;\cbasesb)$ and $\beta \in
    \cfact(A;\cbaseosb)$,
  \end{enumerate}
  and assume that condition \eqref{eq:rtp-iso-cond} holds.  Note that
  by Example \ref{example:gns-cbase}, Lemma \ref{lemma:cbase-gns}, and
  Proposition \ref{proposition:rep-cfact}, given the data listed in
  ii), we can construct the data listed in iii) such that
  \eqref{eq:rtp-iso-cond} is satisfied, and vice versa.

 We can form the classical fiber product $A \rfibrers A$ and
define  representations
\begin{align*}
  \rho_{\leg{2}} &\colon N^{op} \to A \rfibrers A, \ x \mapsto 1 \rtensorm
  \rho(x), & \sigma_{\leg{1}} &\colon N \to A \rfibrers A, \ y
  \mapsto \sigma(y) \rtensorm 1.
\end{align*}
A finite-dimensional {\em Hopf-von Neumann bimodule} is a tuple
$(N,\mu,A,\rho,\sigma,\Delta)$, where $N,\mu,A$, $\rho,\sigma$ are as
above and $\Delta \colon A \to A \rfibrers A$  is a
$*$-homomorphism that satisfies
$\Delta \circ \rho = \rho_{\leg{2}}$ and $\Delta \circ \sigma =
\sigma_{\leg{1}}$ and makes the following diagram commute:
\begin{align*}
  \xymatrix@C=45pt@R=15pt{
    A \ar[r]^{\Delta} \ar[d]^{\Delta} &
    {A \rfibrers A} \ar[d]^{\Delta \underset{\mu}{\ast} \Id} \\
    {A \rfibrers A} \ar[r]^{\Id \underset{\mu}{\ast} \Delta} &
    {A \rfibrers A \rfibrers A.}
  }
\end{align*}

 We can also form the spatial $C^{*}$-fiber product $A \fibreab A$
 and define $C^{*}$-factorizations
\begin{align*}
  \alpha \rt \alpha &:= \lnspan \kalpha{1}\alpha\rnspan \in \cfact(A
  \fibreab A;\cbasesb), &
  \beta \lt \beta &:= \lnspan \kbeta{2} \beta\rnspan \in
  \cfact(A\fibreab A;\cbaseosb),
\end{align*}
where $\kalpha{1}$ and $\kbeta{2}$ were defined below Equation
\eqref{eq:ketbra-ops}. The associated representations are given by
\begin{align*}
 \rho_{(\alpha \rt \alpha)}(b^{\dag}) &= 1 \rtensorh
\rho_{\alpha}(b^{\dag}), & \rho_{(\beta \lt \beta)}(b) &=
\rho_{\beta}(b) \rtensorh 1 
\end{align*}
for all $b^{\dag} \in \frakB$, $b\in \frakB$ \cite[Proposition
2.7]{timmer:cpmu}.  A finite-dimensional {\em concrete Hopf
  $C^{*}$-bimodule} is a tuple $(\cbasesb,H,A,\alpha,\beta,\Delta)$,
where $\cbasesb,H,A,\alpha,\beta$ are as above and $\Delta \in
\Mor(A_{\alpha}, (A \fibreab A)_{\alpha \rt \alpha}) \cap
\Mor(A_{\beta},(A \fibreab A)_{\beta \lt \beta})$ makes the following
diagram commute:
\begin{align*}
  \xymatrix@C=45pt@R=15pt{
    A \ar[r]^{\Delta} \ar[d]^{\Delta} &
    {A \fibreab A} \ar[d]^{\Delta \underset{\frakH}{\ast} \Id} \\
    {A \fibreab A} \ar[r]^{\Id \underset{\frakH}{\ast} \Delta} &
    {A \fibreab A \fibreab A.}
  }  
\end{align*}
Combining the results obtained so far, we find:
\begin{proposition} \label{proposition:hopf-iso}
  Let $\Delta_{\mu} \colon A \to A \rfibrers A$ and $\Delta_{\frakH}
  \colon A \to A \fibreab A$ be
  $*$-homomorphisms such that $\Delta_{\frakH} =
  \phi^{U,\zeta}_{\alpha,\beta} \circ \Delta_{\mu}$. Then
  $(N,\mu,A,\rho,\sigma,\Delta)$ is a Hopf-von Neumann bimodule if and
  only if
 $(\cbasesb,H,A,\alpha,\beta,\Delta)$ is
  a concrete Hopf-$C^{*}$-bimodule. \qed
\end{proposition}
Thus, in the finite-dimensional case, Hopf-$C^{*}$-bimodules and
Hopf-von Neumann bimodules are equivalent descriptions of the same
objects.

\section{Finite-dimensional pseudo-multiplicative unitaries}

In the finite-dimensional case, the notion of a pseudo-multiplicative
unitary and of a $C^{*}$-pseudo-multiplicative unitary
 coincide. To make this statement precise, we recall
the necessary definitions. Let
\begin{enumerate}
\item $H$ be a finite-dimensional Hilbert space,
\item $N$ be a finite-di\-mensional $C^{*}$-algebra with a faithful
  state $\mu$ and nondegenerate faithful representations $\rho \colon
  N^{op} \to {\cal L}(H)$ and $\sigma,\hsigma \colon N \to {\cal
    L}(H)$ such that $\rho(N^{op})$, $\sigma(N)$, $\hsigma(N)$
  commute pairwise,
\item $\cbasesb$ be a finite-dimensional standard $C^{*}$-base with
  bicyclic vector $\zeta \in \frakH$ and $C^{*}$-factorizations $\alpha
  \in \cfact(H;\cbasesb)$, $\beta,\hbeta \in \cfact(H;\cbaseosb)$
  such that $\alpha,\beta,\hbeta$ are pairwise compatible.
\end{enumerate}
 Assume  that
\begin{align} \label{eq:pmu-iso-cond}
    \begin{gathered}
      (\frakH,\frakB,\frakBo) \underset{U}{\sim}
      (H_{\mu},\pi_{\mu}(N),\pi_{\mu}^{op}(N^{op})), \qquad U\zeta =
      \zeta_{\mu},
      \\
      \rho = \rho_{\alpha} \circ \Ad_{U^{*}} \circ \pi_{\mu}^{op}, \quad
      \sigma = \rho_{\beta}  \circ \Ad_{U^{*}} \circ \pi_{\mu}, \quad
      \hsigma = \rho_{\hbeta}  \circ \Ad_{U^{*}} \circ \pi_{\mu}.
  \end{gathered}
\end{align}
Similarly as before, we can, given the data listed in ii),
construct the data listed in iii) such that \eqref{eq:pmu-iso-cond} is
satisfied, and vice versa.

By Proposition \ref{proposition:rtp-iso}, we can identify $\pHsource
\cong \Hsource$ and $\Hrange \cong \pHrange$.  Let
\begin{align*}
  V \colon \pHsource \cong \Hsource \to \Hrange \cong \pHrange
\end{align*}
be a unitary.  Recall that $V$ is a {\em pseudo-multiplicative
  unitary} \cite{vallin:2} if for all $x \in N$, $y \in N^{op}$,
\begin{gather}
  \label{eq:pmu-intertwine}
  \begin{aligned}
    V(\rho(y) \mtimes 1) &= (1 \motimes \rho(y))V, &
    V(\sigma(x) \mtimes 1) &= (\sigma(x) \motimes 1)V, \\
    V(1 \mtimes \sigma(x)) &= (\hsigma(x) \motimes 1)V, & V(1 \mtimes
    \hsigma(x)) &= (1 \motimes \hsigma(x))V,
  \end{aligned}
\end{gather}
and if the following diagram commutes,
\begin{align}  \label{eq:pmu-pentagon}
      \begin{gathered}
        \xymatrix@R=15pt{ {\pHsource \rtensor{\hsigma}{\mu}{\rho} H}
          \ar[r]^{ V \mtimes \Id} \ar[d]_{\Id \mtimes V} & {\pHrange
            \rtensor{\hsigma}{\mu}{\rho} H } \ar[r]^{ \Id \motimes V}&
          { \pHrange \rtensor{\rho}{\mu^{op}}{\sigma} H,}
          \\
          {H \rtensor{\hsigma}{\mu}{\rho_{\leg{2}}} (\pHrange)}
          \ar[d]_{ \Id\mtimes \Sigma^{\mu}} & & { (\pHsource)
            \rtensor{\rho_{\leg{1}}}{\mu^{op}}{\sigma} H} \ar[u]_{ V
            \motimes \Id}
          \\
          {\pHsource \rtensor{\sigma}{\mu}{\rho} H} \ar[rr]^{ V
            \mtimes \Id}&& {(\pHrange)
            \rtensor{\hsigma_{\leg{1}}}{\mu}{\rho} H} \ar[u]_{
            \Sigma^{\mu}_{\leg{23}}} } 
    \end{gathered}
\end{align}
where $\Id  \mtimes \Sigma^{\mu}$ and $\Sigma^{\mu}_{\leg{23}}$ flip the second and
the third component in the respective relative tensor product.

On the other hand, $V$ is a {\em $C^{*}$-pseudo-multiplicative
  unitary} \cite{timmer:cpmu} if it satisfies
 \begin{gather} \label{eq:cpmu-intertwine}
    \begin{aligned}
      V(\alpha \lt \alpha) &= \alpha \rt \alpha, &
      V(\hbeta \rt \beta) &= \hbeta \lt \beta, &
      V(\hbeta \rt \hbeta) &= \alpha \rt \hbeta, & 
      V(\beta \lt \alpha) &= \beta \lt \beta
    \end{aligned}
  \end{gather}
and if the following diagram commutes,
\begin{align}\label{eq:cpmu-pentagon}
      \begin{gathered}
      \xymatrix@R=15pt{ {\Hone} \ar[r]^{ V
          \rtensorh \Id} \ar[d]_{\Id \rtensorh V} & {\Htwo}
        \ar[r]^{ \Id \rtensorh V}& { \Hthree,}
        \\
        {\Hfive} \ar[d]_{ \Id\rtensorh \Sigma^{\frakH}} & & {\Hfour}
        \ar[u]_{ V \rtensorh \Id}
        \\
        {\Hfourlt} \ar[rr]^{ V \rtensorh \Id}&& {\Hfourrt} \ar[u]_{
          \Sigma^{\frakH}_{\leg{23}}} }
    \end{gathered}
\end{align}
where again $\Id \rtensorh \Sigma^{\frakH}$ and
$\Sigma^{\frakH}_{\leg{23}}$ flip the second and the third component
in the respective $C^{*}$-relative tensor product.

Combining the results of the preceding section, we find:
\begin{proposition}
 $V \colon \Hsource \cong \pHsource \to \pHrange \cong \Hrange$ is
  a pseudo-multiplicative unitary if and only if it is a
  $C^{*}$-pseudo-multiplicative unitary. 
  \end{proposition}
  \begin{proof}
    By Proposition \ref{proposition:rep-cfact}, $V$ satisfies
    \eqref{eq:pmu-intertwine} if and only if it satisfies
    \eqref{eq:cpmu-intertwine}. Moreover, using the explicit formula
    for the identifications $\pHsource \cong \Hsource$ and $\Hrange
    \cong \pHrange$ given in Proposition \ref{proposition:rtp-iso},
    one easily verifies that diagram \eqref{eq:pmu-pentagon} commutes
    if and only if diagram \eqref{eq:cpmu-pentagon} commutes.
  \end{proof}
  \begin{remark}
    If the unitary $V$ above is a sufficiently well-behaved
    ($C^{*}$-)pseudo-multiplicative unitary, one can associate to it
    two finite-dimensional Hopf-von Neumann bimodules
    \cite{enock:10} and two finite-dimensional concrete Hopf
    $C^{*}$-bimodules \cite{timmer:cpmu}. One easily verifies that
    these bimodules coincide in the sense of Proposition
    \ref{proposition:hopf-iso}.
 \end{remark}

\def\cprime{$'$}

\end{document}